\documentclass[12pt]{article}
\usepackage{amsfonts}
\usepackage{amsmath}
\usepackage{amsthm}
\newtheorem*{theorem}{Theorem}
\newtheorem*{lemma}{Lemma}
\newtheorem*{corollary}{Corollary}

\theoremstyle{remark}

\newtheorem*{remark}{Remark}

\newtheorem*{comments}{Comments}

\theoremstyle{definition}

\newtheorem*{definition}{Definition}

\def\rom#1{{\rm#1}}

\def\Z{{{\mathbb Z}}}
\def\C{{\mathbb C}}

\def\gl{{\mathfrak{gl}}}
\def\sl{{\mathfrak{sl}}}
\def\n{{\mathfrak{n}}}
\def\h{{\mathfrak{h}}}
\def\a{{\mathfrak{a}}}
\def\hgl{\widehat{\gl}_\infty}
\def\gli#1{\gl_{\frac\infty2}^{(#1)}}
\def\Det{\text{\rm Det}}
\def\Ind{\text{\rm Ind}}
\def\oplusl{\mathop\oplus\limits}
\def\bopl{\mathop\bigoplus\limits}
\def\botiml{\mathop\bigotimes\limits}

\def\wA{\widetilde A}
\def\wD{\widetilde D}

\def\CG{Clebsch -- Gordan\ }

\newcounter{tmp}
\newenvironment{enums}{\begin{list}{$(\arabic{tmp})$}%
{\usecounter{tmp}}}{\end{list}}
\newenvironment{lgather*}{\begin{gather*}}{\end{gather*}}

\begin{document}

\author{by BORIS SHOIKHET, e-mail: burman@ium.ips.ras.ru}
\title{Representations   of   the   Lie   algebra   $\hgl$   and
``reciprocity formula'' for \CG coefficients}
\date{}

\maketitle

\section*{Introduction}

In this article I compare two approaches to the problem of
calculating the characters of the Lie algebra~$\hgl$ irreducible
representations with ``semidominant'' highest weights
and integral central charge. The first approach is the
remarkable result of V.~Kac and A.~Radul~\cite1, the second one
is a small generalization of the author's approach in~\cite2. As
central charge $c\in\Z_{\le0}$ tends to~$-\infty$, the second
approach also gives a precise answer. In this case, both
answers represent $\hgl$-module as a direct sum of the irreducible
$\gli1\oplus\gli2$-modules
($\gli i$ are natural subalgebras
in~$\hgl$), but in completely different form. The equivalence
of both formulas gives some ``reciprocity formula'' for
\CG coefficients (see~\eqref{eq6} of~\ref{s2}). These
coefficients are defined by the identity $L(\alpha)\otimes
L(\beta)=\oplusl_\gamma C_{\alpha\beta}^\gamma L(\gamma)$,
$C_{\alpha\beta}^\gamma\in \Z_{\ge0}$, for tensor product of two
irreducible~$\gl_N$(or~$\gl_{\frac\infty2}$)-modules with dominant highest
weights.
(These coefficients depend of~$N$, but they stabilize when~$N$
tends to~$\infty$).

I am grateful to S.~Arkhipov for useful discussions.

\section{Decomposition of the induced representations}
\label{s1}

\subsection{}\label{1.1}
There are two natural subalgebras, $\gl^{(1)}_n$ and
$\gl^{(2)}_n$, in the Lie algebra $\gl_{2n}$ of complex
$2n\times2n$-matrices, and there are also two Abelian
subalgebras, $\a_+$ and $\a_-$ (see Fig.~\ref{f1}).

\begin{figure}[h]
$$
\vcenter{\hbox{$
\addtolength\arraycolsep{5mm}
\renewcommand\arraystretch{2.5}
\begin{array}{|c|c|}
\hline
\gl_n^{(1)}&\a_+\\
\hline
\a_-&\gl_n^{(2)}\\
\hline
\end{array}$}}\quad
2n
$$
\caption{}\label{f1}
\end{figure}

Adjoint action of the subalgebra $\gl_n^{(1)}\oplus\gl_n^{(2)}$
preserves the universal enveloping algebra $U(\a_-)\cong
S^*(\a_-)$, and we want to decompose~$S^*(\a_-)$ with respect to
this action. Let $\n_-\oplus\h\oplus\n_+=\gl_{2n}$ and
$\n_-^{(i)}\oplus\h^{(i)}\oplus\n_+^{(i)}=\gl_n^{(i)}$ be the
standart Cartan decompositions, and let $E_{ij}$ be the element
of~$\gl_{2n}$ with~$1$ in the $(i,j)$-cell and $0$ in other cells. Next,
denote
$$
\Det_k=\Det
\begin{pmatrix}
E_{n+1,n-k+1}&\hdots&E_{n+1,n}\\
\vdots&&\vdots\\
E_{n+k,n-k+1}&\hdots&E_{n+k,n}
\end{pmatrix}
\in U(\a_-)\cong S^*(\a_-)\quad(k=1\ldots n)
$$

\subsubsection{}\label{1.1.1}
\begin{lemma}
Monomials $\Det^{l_1}_1\cdot\ldots\cdot\Det^{l_n}_n$
$(l_i\in\Z_{\ge0})$ exhaust all the
$\gl_n^{(1)}\oplus\gl_n^{(2)}$-singular vectors with respect to
the adjoint action.
\end{lemma}

\subsubsection{}\label{1.1.2}
\begin{lemma}
With respect to the adjoin action of the Lie algebra
$\gl_n^{(1)}\oplus\gl_n^{(2)}$, $S^*(\a_-)$ decomposes into
direct sum of finite-dimensional irreducible modules $L_w\otimes
L_w$ where $w$ goes through all the monomials
$\Det^{l_1}_1\cdot\ldots\cdot\Det^{l_n}_n$,
$l_i\in\Z_{\ge0}$.
\end{lemma}

We will prove these Lemmas in Section~\ref{1.2}.

\subsubsection{}\label{1.1.3}
\begin{remark}

In the sequel we will need an expression of the highest weight
$\theta_{l_1,\ldots,l_n}$ of monomial
$\Det^{l_1}_1\cdot\ldots\cdot\Det^{l_n}_n$ through the basis
$\{E_{ii}\}$ of the Cartan subalgebras $\h^{(1)}\oplus\h^{(2)}$.
We have:
\setbox57=\hbox{$\theta_{l_1,\ldots,l_n}(E_{22})= -l_n-l_{n-1}$}
$$
\begin{aligned}
{}&\theta_{l_1,\ldots,l_n}(E_{11})=-l_n\\
{}&\theta_{l_1,\ldots,l_n}(E_{22})= -l_n-l_{n-1}\\
\noalign{\hbox to \wd57 {\dotfill}}
{}&\theta_{l_1,\ldots,l_n}(E_{nn})=-l_n-l_{n-1}-\ldots-l_1
\end{aligned}
$$
and
$$
\begin{aligned}
{}&\theta_{l_1,\ldots,l_n}(E_{n+1,n+1})=l_1+\ldots+l_n\\
\noalign{\hbox to \wd57 {\dotfill}}
{}&\theta_{l_1,\ldots,l_n}(E_{2n,2n})=l_n
\end{aligned}
$$
The corresponding Young diagram is shown in the  Fig.~\ref{f2}:
\begin{figure}[h]
\begin{center}
\unitlength=1em
\begin{picture}(23,10)
\put(1,3.5){$n$}
\put(8.2,4){$n-1$}
\put(4.5,8.3){$l_n$}
\put(9,8.3){$l_{n-1}$}
\put(19,8.3){$l_1$}
\put(14,6.5){$\dots$}
\put(2,0){\line(0,1){8}}
\put(8,0){\line(0,1){8}}
\put(12,1){\line(0,1){7}}
\put(15,2){\line(0,1){3}}
\put(17,5){\line(0,1){3}}
\put(22,7){\line(0,1){1}}
\put(2,0){\line(1,0){6}}
\put(8,1){\line(1,0){4}}
\put(12,2){\line(1,0){3}}
\put(15,5){\line(1,0){2}}
\put(17,7){\line(1,0){5}}
\put(2,8){\line(1,0){20}}
\end{picture}
\end{center}
\vskip1em
\caption{}\label{f2}
\end{figure}

\end{remark}

\subsection{}\label{1.2}
In this Section we prove Lemmas \ref{1.1.1} and \ref{1.1.2}.

\subsubsection{}
\label{1.2.1}
\begin{proof}[Proof of Lemma \ref{1.1.1}\/\rom:]

The fact that the monomials
$\Det^{l_1}_1\cdot\ldots\cdot\Det^{l_n}_n$ are
$\gl_n^{(1)}\oplus\gl_n^{(2)}$-singular vectors, is
straightforward. Conversely, let $\xi$ be a
$\gl_n^{(1)}\oplus\gl^{(2)}_n$-singular vector, $\xi\in
S^*(\a_-)$. Consider the minimal rectangular domain in~$\a_-$
with a vertex in the central (upper right) corner which includes
all the elements contained in the notation of~$\xi$. Suppose,
for example, that its horizontal side is not smaller then the
vertical one, and has the length~$l$. Then there are no more
than~$l$ squares in the $l$th column, and elements of $\n_+^{(i)}$
shift the $l$th column to $1$st, $\ldots$, $(l-1)$th. Applying
these elements, we should obtain~$0$, and we obtain ($l-1$)
equations on~$\xi$. Therefore, if $\xi$ is {\it linear\/} in
elements of the $l$th column, then $\xi=C\cdot\Det_l$,
where $C$ is an
expression in a smaller square.
Next, all (not
only linear) expressions in elements of $l$th column that vanish
under the corresponding $(l-1)$ vector fields, are divisible on
some degree of $\Det_l$, say $\Det_l^k$, via arguments of
grading, and we have $\xi=C\cdot\Det^k_l$,
where $C$ is an
expression in a smaller square. We can apply the preceding
arguments to the expression~$C$. (see also~\cite2, Ch.~I, \S1,
Sec.~1.3)
\end{proof}

\subsubsection{}\label{1.2.2}
\begin{proof}[Proof of Lemma~\ref{1.1.2}]

The $\n_-^{(1)}\oplus\n_-^{(2)}$-action on $S^*(\a_-)$ is locally nilpotent.
\end{proof}

\subsubsection{}\label{1.2.3}
In the case $n=\infty$ we obtain the following result:
\begin{equation}\label{eq1}
S^*(\a_-)
=\bopl_%
{\genfrac{}{}{0pt}{}{{\mathrm{for\ all}}}{{\mathrm{Young\ diagrams\ }}D}}
L_D\otimes L_D
\end{equation}
where $L_D$ is the irreducible $\gl_\frac\infty2$-module
with highest weight~$\theta_D$.

\subsection{}\label{1.3}

\begin{definition}
The weight $\chi\colon\h\to\C$ is a {\it semidominant}, if its restrictions
$\chi^{(1)}$ and $\chi^{(2)}$ on subalgebras $\gl_n^{(1)}$ and
$\gl_n^{(2)}$ are dominant weights. In other words, in the basis
$\{E_{ii}\}$ of~$\h$ we have $\chi=(\chi_1,\ldots,\chi_n)$ where
all the~$\chi_i$ belong to~$\Z$,
$\chi_1\ge\chi_2\ge\ldots\ge\chi_n$, and
$\chi_{n+1}\ge\ldots\ge\chi_{2n}$.
\end{definition}

Difference
$\chi_n-\chi_{n+1}$ may be a negative integer.

\subsubsection{}\label{1.3.1}
Let $\chi\colon\h\to\C$ be a semidominant weight, and let
$L_{\chi^{(1)}}\otimes L_{\chi^{(2)}}$ be corresponding (finite
dimensional) irreducible $\gl_n^{(1)}\oplus\gl_n^{(2)}$-module.
We can continue this module on subalgebra
$\gamma=\h+\a_+\oplus\gl_n^{(1)}\oplus\gl_n^{(2)}$ by formulas
$hv=\chi(h)v$ for any~$h\in\h$, and $\a_+v=0$.

\begin{definition}
$\Ind_\chi=U(\gl_{2n})\botiml_{U(\gamma)}(L_{\chi^{(1)}}\otimes L_{\chi^{(2)}})$.
\end{definition}

\subsubsection{}\label{1.3.2}
\begin{lemma}
As a $\gl_n^{(1)}\oplus\gl_n^{(2)}$-module, $\Ind_\chi$ is
$(L_{\chi^{(1)}}\otimes L_{\chi^{(2)}})\bigotimes S^*(\a_-)$.
\end{lemma}

In particular, $\Ind_\chi$ decomposes into the direct sum of
finite-dimensional irreducible
$\gl_n^{(1)}\oplus\gl_n^{(2)}$-modules.

\begin{proof}[Proof\/\rom:] It is obvious.
\end{proof}

\subsubsection{}\label{1.3.3}
Consider the case of~$\hgl$.

With all the values $c$ of central charge and $\chi_{centr}$ of
highest weight~$\chi$ on ``central'' coroot, as
$\gli1\oplus
\gli2$-module.
\begin{multline}\label{eq2}
\Ind_{\chi,c}
=\bopl_%
{\genfrac{}{}{0pt}{}{{\mathrm{for\ all}}}{{\mathrm{Young\ diagrams\ }}D}}
(L_{\chi^{(1)}}\otimes L_{\chi^{(2)}})\bigotimes(L_D\otimes L_D)=\\
=\bopl_D
(L_{\chi^{(1)}}\otimes L_D)\bigotimes
(L_{\chi^{(2)}}\otimes L_D)
\end{multline}
We have:
\begin{equation}\label{eq3}
L_\alpha\otimes L_\beta=\oplusl_\gamma C_{\alpha\beta}^\gamma L_\gamma,
\end{equation}
where $L_\alpha$, $L_\beta$, $L_\gamma$ are irreducible
$\gl_\frac\infty2$-modules with highest weights $\alpha$,
$\beta$,~$\gamma$; $C^\gamma_{\alpha\beta}\in\Z_{\ge0}$ are
called \CG coefficients.

Therefore, as $\gli1\oplus\gli2$-module,
\begin{equation}\label{eq4}
\Ind_{\chi,c} =\bopl_{D,\nu_1,\nu_2}C_{\chi^{(1)},D}^{\nu_1}
C_{\chi^{(2)},D}^{\nu_2}L_{\nu_1}\otimes L_{\nu_2}
\end{equation}
In particular,  the  multiplicity  of  $\gli1\oplus\gli2$-module
$L_{\nu_1}\otimes L_{\nu_2}$ in $\Ind_{\chi,c}$ is equal to
\begin{equation}\label{eq5}
\sum_{\text{for all $D$}}
C_{\chi^{(1)},D}^{\nu_1}\cdot
C_{\chi^{(2)},D}^{\nu_2}
\end{equation}
Note, that the sum is finite.

\section{Irreducible  representations of the
Lie algebra $\hgl$ with
$c=-N$ when $N\to\infty$}
\label{s2}

\subsection{}\label{2.1}
In \ref{s1} we decomposed the  $\hgl$-module  $\Ind_{\chi,-N}$
into the direct  sum of irreducible $\gli1\oplus\gli2$-modules, and
this decomposition does not depend  on  central  charge~$c$  and
$\chi_{centr}=\chi(\alpha_0^\vee)$;  and  it  is  interesting to
find irreducible factor of $\Ind_{\chi,-N}$ in the same  terms,
i.e. ``point out'' $\gli1\oplus\gli2$-modules, which are
present in this irreducible factor, with dependence on $c$  and
$\chi_{centr}$.  This  is  very  interesting, and, I think, very
difficult problem. It was solved in  the  author's  paper~\cite2
only  in  a  particular  case,  when   $\chi=0$   and any $c\in\C$
(see~\ref{2.1.5} and \cite2, Ch.~I, \S2).

However, when $c=-N\in\Z_{\le0}$ and $N\to\infty$, corresponding
irreducible factor of the
module $\Ind_{\chi,-N}$ tends to the whole
module $\Ind_{\chi,-N}$. Compairing this fact with the result of
V.~Kac   and   A.~Radul~\cite1   we   obtain  relations  on  \CG
coefficients.

\subsubsection{}\label{2.1.1}
\begin{theorem}
Let central charge $c=-N$  $(N\gg0)$.  Then  there  exists
$\delta(N)\in\Z_{\ge0}$  which  tends  to~$\infty$  as $N$ tends
to~$\infty$, and such that there are no $\hgl$-singular  vectors
on  all the levels of the representation $\Ind_{\chi,-N}$, which
are less then $\delta(N)$.
\end{theorem}

\begin{corollary}
Maximal submodule in $\Ind_{\chi,-N}$ does  not  intersect  with
levels less than~$\delta(N)$.\qed
\end{corollary}

Only  $\gli1\oplus\gli2$-singular  vector may be $\hgl$-singular
vector,   this   vector   is   equal   to   the   sum   of   the
$\gli1\oplus\gli2$-singular vectors in several representations
$(L_{\chi^{(1)}}\otimes L_D)\bigotimes
(L_{\chi^{(2)}}\otimes   L_D)$.   We   want   to  show,  that
diagram~$D$ should be very big when~$N\gg0$.

\subsubsection{}\label{2.1.2}
\begin{lemma}  Let  $V$,  $W$  be  two  irreducible   $\gl_N$\rom(or
$\gl_{\frac\infty2}$\/\rom)-modules,  let $v$ and $w$ be their highest vectors.
Then             every              singular              vector
$\sum\limits_i\theta_1^{(i)}v\otimes\theta_2^{(i)}v$    contains
the term with~$\theta_1=1$.
\end{lemma}

\begin{proof}[Proof\/\rom:]    $V$    and    $W$    are    irreducible
representations, and hence they do not contain singular vectors
besides their highest weight vectors.
\end{proof}

\subsubsection{}\label{2.1.3}
Let  $e_0\in\n_+\subset\gl_{2n}(\hgl)$ be the ``corner''
element   of~$\a_+$   (see   Fig.~\ref{f1}),   in    the    case
of~$\gl_{2n}$, $e_0=E_{n,n+1}$; and let $\alpha_0^\vee$ be
``central''   coroot,  $\chi(\alpha_0^\vee)=\chi_{centr}$.  We
will need in the sequel the direct expression for the commutator
$[e_0,\Det_k^l]$. For this, introduce the  following  notations.
Denote  by  $\{y_{ij}\}$  the  elements  of~$\a_-$,  in the case
of~$\gl_{2n}$  $y_{ij}=E_{n+i,n-j+1}$;   we   will   use   these
notations     also     in    the    case    of~$\hgl$.    Denote
$A_k=(y_{ij})_{i,j=1\ldots k}$, $\wA_k=(y_{ij})_{i,j=2\ldots k}$
(we   have   $\Det_k=\Det   A_k$).   Let   $\wA_{ij}$   be   the
matrix~$\wA_k$ without $j$th row and $i$th column. Let $z_i^+$
be  the element from~$\n_-^{(1)}$, standing in the $i$th column
above~$y_{1i}$, and  $z_j^-$  be  an  element  from~$\n_-^{(2)}$,
standing  in the $j$th row
right from~$y_{j1}$. Now we are ready to
formulate the result:

\begin{lemma}
\begin{multline*}
[e_0,\Det_k^l]=(-1)^{1+k}\cdot l\cdot(\Det\wA_k\cdot\Det_k^{l-1}
(\alpha_0^\vee+c+k-l)+\\
+\sum_{i,j=2\ldots k}(-1)^{i+j}\cdot l\cdot\Det\wA_{ij}\cdot\Det_k^{l-1}
(y_{j1}z_i^+-y_{1i}z_j^-)
\end{multline*}
\end{lemma}

\begin{proof}[Proof.]
It is a direct calculation from~\cite2, Ch.~I, \S1.
\end{proof}

Note, that the sum belongs to
$U(\a_-)\n_-^{(1)}\oplus
U(\a_-)\n_-^{(2)}$.

\subsubsection{}\label{2.1.4}
\begin{proof}[The proof of the Theorem\/\rom:]

Any $\hgl$-singular vector is represented as
$$
\theta=\sum_D\sum_i[n_{s_i}^{i,D}[\cdots[n_1^{i,D},D]\cdots]\cdot
n_{s_{i+1}}^{i,D}\cdot\ldots\cdot n_{r_i}^{i,D}\cdot v,
$$
where $n_j^{i,D}\in\n_-^{(1)}\oplus\n_-^{(2)}$ and  $v$  is  the
highest  weight  vector.  We want to prove that $[e_0,\theta]\ne0$
for $N\gg0$.

\begin{enums}
\item
$e_0$ commutes with $\n_-^{(1)}\oplus\n_-^{(2)}$, and so,
$$
[e_0,\theta]=\sum_D\sum_i[n_{s_i}^{i,D}[\cdots
[n_1^{i,D}[e_o,D]\cdots]
n_{s_{i+1}}^{i,D}\cdot\ldots\cdot n_{r_i}^{i,D}\cdot v
$$
\item
by Lemma~\ref{2.1.2}, the  notation  of~$\theta$  contains
the      term     $D\cdot\alpha     v$,     where     $\alpha\in
U(\n_-^{(1)}\oplus\n_-^{(2)})$, we want to prove, that the  term
$[e_0,D]\alpha\cdot v$ can not annihilate with other terms
\item
We will consider only the case $D=\Det_k^l$, the general case is similar
\item
Denote $\wD=\Det\wA_k\Det_k^{l-1}$; by Lemma~\ref{2.1.3},
$$
[e_0,\Det_k^l]=(-1)^{1+k}\cdot l\cdot\wD(\alpha_0^\vee+c+k-l)+o(N),
$$
the first term is very big, because $c$ acts as multiplication by~$-N$
\item
For $n\in\n_-^{(1)}\oplus\n_-^{(2)}$ we have:
$$
[n,\wD(\alpha_0^\vee+c+k-l)]=[n,\wD](\alpha_0^\vee+c+k-l)+\wD[n,\alpha_0^\vee]
;
$$
the second summand is $o(N)$, and $[n,\wD]$ is not equal to~$\wD$
element of~$U(\a_-)$
\end{enums}
\end{proof}

\subsubsection{}\label{2.1.5}
\begin{remark}
It is possible to get the explicit answer for  the  $\hgl$-module
$\Ind_\mu$    (with    $\chi=0$    and    $c=\mu\in\C$):    when
$\mu=N\in\Z_{\ge0}$, the singular vectors  are  $\Det_1^{N+1}v$,
$\Det_2^{N+2}v$,  $\ldots$,  and  when $\mu=-N\in\Z_{\le0}$, the
singular vectors are $\Det_{N+1}v$, $\Det^2_{N+2}v$, $\ldots$. We
see that irreducible factor of $\Ind_\mu$  ($\mu=\pm N$)
tends  to  $\Ind_\mu$
when  $N\to\infty$.  When  $\mu\not\in\Z$  there are no singular
vectors in $\Ind_\mu$, and $\Ind_\mu$ is irreducible. This is  a
result from~\cite2, Ch.~I,~\S1.
\end{remark}

\subsection{}\label{2.2}
In  this Section, we compare Theorem~\ref{2.1.1} with result from
the
paper by V.~Kac and A.~Radul~\cite1.

\subsubsection{}\label{2.2.1}
We define a set
$$
H_N=\{(\nu_1,\nu_2,\ldots,\nu_N),\,\nu_1\ge\nu_2\ge\ldots\ge\nu_N,\,
\nu_i\in\Z\quad\text{for all}\quad i=1\ldots N\},
$$
$H^+_N$ ($H^-_N$) is
a  subset  in  $H_n$  which  consists   of   $\nu_i\in\Z_{\ge0}$
($\nu_i\in\Z_{\le0}$).  Let  $\nu\in  H_N$; define $\hgl$-weight
$\Lambda(\nu)$                    as                    follows:
$$
\Lambda(\nu)=(\ldots0,0,\nu_{p+1},\ldots,\nu_N;
\nu_1,\nu_2,\ldots,\nu_p,0,0,\ldots),
$$
where
$\nu_1,\ldots,\nu_p\ge0$, $\nu_{p+1},\ldots,\nu_N\le0$,  and  we
put a semicolon between the 0th and the first slot. So,
$\Lambda_-(\nu)=(\ldots0,0,\nu_{p+1},\ldots,\nu_N)$ is the $\gli1$-weight
and
$\Lambda_+(\nu)=(\nu_1,\nu_2,\ldots,\nu_p,0,0,\ldots)$ is the $\gli2$-weight.
\begin{theorem}[\cite1]
As $\gli1\oplus\gli2$-module,
$$
L(\Lambda(\nu),-N)=\bopl_{\genfrac{}{}{0pt}{}{\lambda\in H^-_N}{\mu\in H^+_N}}
C^\nu_{\lambda,\mu}
L_-(\lambda)\otimes L_+(\mu).
$$
\end{theorem}

\begin{comments}
$L(\Lambda(\nu),-N)$  is the irreducible $\hgl$-module with the
highest weight $\Lambda(\nu)$ and the central charge~$-N$; where
$C_{\lambda\mu}^\nu$  are   \CG   coefficients   for   the   Lie
algebra~$\gl_N$  (see~\eqref{eq3});  $L_-(\lambda)(L_+(\mu))$ are
the irreducible $\gli1$($\gli2$)-modules with highest weights~$\lambda$
(resp.~$\mu$).
\end{comments}

\subsubsection{}\label{2.2.2}
Send $N$ to~$\infty$, and choose a $\gl_N$-weight~$\nu$ such that
$$
\nu=(\nu_1,\ldots,\nu_k,0,0\ldots\ldots,0,0,\nu_s,\ldots,\nu_N)
$$
where  $\nu_1,\ldots,\nu_k>0$,  $\nu_s,\ldots,\nu_N<0$;  as  $N$
grows we just add zeros between them.

\begin{theorem}
\begin{multline}\label{eq6}
\sum_%
{\genfrac{}{}{0pt}{}{{\mathrm{for\ all}}}{{\mathrm{Young\ diagrams\ }}D}}
C^{\lambda_-}_{\Lambda_-(\nu),D_-}\cdot
C^{\mu_+}_{\Lambda_+(\nu),D_+}\quad(\text{in the sense of}\
\gl_{\frac\infty2})\,=\\
=C^\nu_{\lambda_-,\mu_+}\quad(\text{in the sense of $\gl_N$,\
$N\gg0$, $\lambda_-$, $\mu_+$, and  $\nu$ are fixed})
\end{multline}
\end{theorem}

\begin{comments}
If $D$ is as in Fig.~\ref{f1}, then
$$
D_-=\theta_{D_-}=
(\ldots,0,0,-l_n,\ldots,-l_2-\ldots-l_n,-l_1-\ldots-l_n)
$$
and
$$
D_+=\theta_{D_+}=
(l_1+\ldots+l_n,l_2+\ldots+l_n,\ldots,l_n,0,0,\ldots);
$$
$\Lambda_+(\nu)$   and   $\Lambda_-(\nu)$   were   defined    in
\ref{2.2.1}.  The  sum  in  the  left-hand  side of~\eqref{eq6} is
finite. The right-hand side stabilizes, as $N$ grows.
\end{comments}

\begin{proof}[Proof\/\rom:]
it is an obvious consequence of Theorems \ref{2.1.1} and \ref{2.2.1}
and~\eqref{eq5} of~\ref{s1}.
\end{proof}

\subsubsection{}\label{2.2.3}
\begin{remark}
For the tautological $\sl_N$-module~$V$,
 tensor  products  $V\otimes  V$  and $V\otimes V^*$ are not
isomorphic:
\begin{align*}
{}&V\otimes V=S^2(V)\oplus\Lambda^2(V),\\
{}&V\otimes V^*=\sl_N\oplus\C
\end{align*}
The   $\gl_N$-weights    $\lambda_-$ (resp.~$\mu_+$)   define
representations  in  tensor  powers  of~$V^*$  (resp.~$V$),  and
weight~$\nu$ defines ``mixed'' representation.
\end{remark}

\end{document}